\documentclass[12pt, reqno]{amsart}
\usepackage{amsmath, amsthm, amscd, amsfonts, amssymb, graphicx, color}
\usepackage[bookmarksnumbered, colorlinks, plainpages]{hyperref}
\hypersetup{colorlinks=true,linkcolor=red, anchorcolor=green, citecolor=cyan, urlcolor=red, filecolor=magenta, pdftoolbar=true}

\newtheorem{thm}{Theorem}[section]
\newtheorem{lem}[thm]{Lemma}

\newtheorem{prop}[thm]{Proposition}
\newtheorem{cor}[thm]{Corollary}
\newtheorem{defn}[thm]{Definition}
\newtheorem{rem}[thm]{Remark}
\newtheorem{rem-ntn}[thm]{Remark and Notation}

\newenvironment{prf}{{\noindent \textbf{Proof:}\ }}{\hfill $\Box$\\ \smallskip}

\numberwithin{equation}{section}

\newcommand{\smnoind}{{\smallskip\noindent}}

\newcommand{\id}{{\rm id}}

\newcommand{\norm}[1]{\left\|#1\right\|}

\newcommand{\botimes}{\bar{\otimes}}

\newcommand{\CL}{\mathcal{L}}
\newcommand{\CK}{\mathcal{K}}

\newcommand{\CA}{A}

\newcommand{\KH}{\mathfrak{H}}

\newcommand{\KK}{\mathfrak{K}}

\newcommand{\BC}{\mathbb{C}}

\newcommand{\BG}{{\mathbb{G}}}

\newcommand{\Bo}{\mathbf{1}}

\newcommand{\BGd}{{\widehat{\mathbb{G}}}}

\newcommand{\Rep}{\mathrm{Rep}}

\newcommand{\ru}{\mathrm{u}}

\newcommand{\crp}{\mathrm{Corep}}

\newcommand{\intw}{\mathrm{Intw}}

\newcommand{\circledT}{\raisebox{.6pt}{\textcircled{\raisebox{-.5pt} {$\scriptstyle\top$}}}}

\begin{document}

\setcounter{page}{1}

\title[Bekka-type amenabilities for unitary corepresentations]{Bekka-type amenabilities for unitary corepresentations of locally compact quantum groups}

\author[X. Chen]{Xiao Chen}

\address{School of Mathematics and Statistics, Shandong University, Weihai, Shandong Province 264209, China.}
\email{\textcolor[rgb]{0.00,0.00,0.84}{cxwhsdu@126.com; chenxiao@sdu.edu.cn}}

%\address{$^{2}$Department of Pure Mathematics, Ferdowsi University of Mashhad, P. O. Box 1159, Mashhad 91775, Iran;
%\newline
%Tusi Mathematical Research Group (TMRG), Mashhad, Iran.}
%\email{\textcolor[rgb]{0.00,0.00,0.84}{second@afa.ac.ir}}

%\dedicatory{This paper is dedicated to Professor ABCD}

\let\thefootnote\relax\footnote{Copyright 2016 by the Tusi Mathematical Research Group.}



\date{Received: 24 January 2017; Revised: 10 April 2017; Accepted: 4 May 2017.
%\newline \indent $^{*}$Corresponding author
}

%-----------------------------------------------------------------------------------------

\begin{abstract}
In this short note, further to Ng's study, we extend Bekka amenability and weak Bekka amenability to general locally compact quantum groups.
We generalize some Ng's results to the general case.
In particular, we show that, a locally compact quantum group $\BG$ is co-amenable if and only if the contra-corepresentation of its fundamental multiplicative unitary $W_\BG$ is Bekka amenable,
and $\BG$ is amenable if and only if its dual quantum group's fundamental multiplicative unitary $W_\BGd$ is weakly Bekka amenable.
\end{abstract} \maketitle

%-----------------------------------------------------------------------------------------

\section{Introduction}

The notion of amenability essentially begins with Lebesgue (1904). In 1929, von Neumann introduced and studied the class of amenable groups and used it to explain why the Banach-Tarski Paradox occurs only for dimension greater than or equal to three.
In 1950, Dixmier extended the concept of amenability to topological groups (see \cite{Pie} and \cite{Vok}).
In 1970s, amenability and co-amenability for Kac algebras were introduced by D. Voiculescu, studied further by M. Enock and J.-M. Schwartz and later by Z.-J. Ruan (see \cite{Voi} and \cite{Ruan}).
In \cite{Ng}, following Bekka's paper \cite{Bek}, C.-K. Ng introduced Bekka amenability and weak Bekka amenability for unitary co-representations of Kac algebras, and used them to characterize amenability and co-amenability for Kac algebras.
Later, amenability and co-amenability for Hopf $C^*$-algebras was investigated by C.-K. Ng (see \cite{Ng00} and \cite{Ng02}).
In 2003, E. B$\acute{e}$dos, R. Conti and L. Tuset extended amenability and co-amenability to algebraic quantum groups and locally compact quantum groups (see \cite{BCT} and \cite{BT}).

%\smallskip

In this short note, we give some remarks on Ng's paper \cite{Ng}. We extend Bekka amenability and weak Bekka amenability to general locally compact quantum groups.
Furthermore, we prove that a locally compact quantum group $\BG$ is co-amenable if and only if the contra-corepresentation of its fundamental multiplicative unitary $W_\BG$ is Bekka amenable,
and $\BG$ is amenable if and only if its dual group's fundamental multiplicative unitary $W_\BGd$ is weakly Bekka amenable.
These results generalize the corresponding propositions for Kac algebras in Ng's paper \cite{Ng}.

%\smallskip

The notions of Bekka-type amenabilities, studied in this note, originate from Bekka's paper \cite{Bek}.
In the case of locally compact groups, all of Bekka-type amenabilities for unitary corepresentations are equal to amenability (introduced by Bekka in \cite{Bek}) for unitary representations.
Remarkably, Bekka showed, in \cite{Bek}, that amenability for a locally compact group is equivalent to the fact that every unitary representation is amenable.
These justify the use of the term ``Bekka-type amenabilities".

%\smallskip

This note is organized as follows. After some preliminaries in Section 2, we discuss  in Section 3 Bekka amenability and weak Bekka amenability for locally compact quantum groups.

%\bigskip

\section{Notations and definitions}

\subsection{Some notations}\label{subsec:ntn}

In this note, we use the convention that the inner product $\langle \cdot , \cdot \rangle$ of a complex Hilbert space $\KH$ is conjugate-linear in the second variable.
We denote by $\CL(\KH)$ and $\CK(\KH)$ the set of bounded linear operators and that of compact operators on $\KH$, respectively.
For any $x,y\in\KH$ and $T\in\CL(\KH)$, we denote by $\omega_{x,y}$ the normal functional given by $$\omega_{x,y}(T):=\langle Tx,y\rangle.$$

The symbol $\otimes$ denotes either a minimal $C^*$-algebraic tensor product or a tensor product of Hilbert spaces, and $\botimes$ denotes a von Neumann algebraic tensor product. 
Moreover, we denote by $\id$ the identity map. Finally, if $X$ and $Y$ are $C^*$-algebras or Hilbert spaces, we use the symbol $\Sigma$ to denote the canonical flip map from $X\otimes Y$ to $Y\otimes X$ sending $x\otimes y$ onto $y\otimes x$, for all $x\in X$ and $y\in Y$. Note that $\Sigma^2=\id$.

For a $C^*$-algebra $\CA$, we use $\Rep(\CA)$ to denote the collection of unitary equivalence classes of non-degenerate $*$-representations of $\CA$.
Let us also recall some notations concerning $\Rep(\CA)$.

Suppose that $(\mu, \KH),(\nu, \KK)\in \Rep(A)$.
We write $\nu\prec \mu$ if $\ker \mu\subset \ker \nu$.

\subsection{Locally compact quantum group}

Let $(C_0(\BG), \Delta, \varphi,\psi)$ be a reduced locally compact quantum group as introduced in \cite[Definition 4.1]{KV1} (for simplicity, we denote it by $\BG$).
The dual locally compact quantum group of $\BG$ is denoted by $(C_0(\BGd), \widehat{\Delta}, \widehat{\varphi},\widehat{\psi})$ (or simply, $\BGd$).
We use $L^2(\BG)$ to denote the Hilbert space given by the GNS construction of the left invariant Haar weight $\varphi$ and consider both $C_0(\BG)$ and $C_0(\widehat{\BG})$ as $C^*$-subalgebras of $\CL(L^2(\BG))$. Notice that $L^2(\BG)=L^2(\BGd)$.

%\smallskip

Let $\Bo$ be the identity of $M(C_0(\BG))$.
There is a unitary $$W_\BG\in M(C_0(\BG)\otimes C_0(\BGd))\subseteq \CL(L^2(\BG)\otimes L^2(\BG)),$$ called the \emph{fundamental multiplicative unitary}, that implements the comultiplication: $$\Delta(x)=W_\BG^*(\Bo\otimes x)W_\BG \quad (x\in C_0(\BG)).$$
We denote by $W_\BGd$ the fundamental multiplicative unitary for the dual quantum group $\BGd$ given by $\Sigma W^*_\BG\Sigma$, where $\Sigma$ is the flip map as defined above. For more details, the readers may refer to \cite{KV1} and \cite{Tim}.

%\smallskip

The von Neumann subalgebra $L^{\infty}(\BG)$ generated by $C_0(\BG)$ in $\CL(L^2(\BG))$ is a Hopf von Neumann algebra under a comultiplication $\widetilde{\Delta}$ defined by $W_\BG$ as in the above (see \cite{KV2} or \cite[Section 8.3.4]{Tim}).
We usually call $L^{\infty}(\BG)$ the von Neumann algebraic quantum group of $\BG$.
Then $L^1(\BG)$ denotes the predual of $L^{\infty}(\BG)$, and $L^1_*(\BG):=\{\omega\in L^1(\BG)\ |\ \exists\eta\in L^1(\BG)\ s.t.\ (\omega\otimes\id)(W_{\BG})^*=(\eta\otimes\id)(W_{\BG}) \}$ is a dense $^*$-subalgebra of $L^1(\BG)$ as introduced in \cite[Page 294-295]{Kus}.

%\medskip

\subsection{Corepresentation}

For any Hilbert space $\KH_U$, a unitary $U\in M(\CK(\KH_U)\otimes C_0(\BG))$ is called a \emph{unitary corepresentation} of $\BG$ on $\KH_U$ if
\begin{equation}\label{eqt:defn-corep}
(\id\otimes\Delta)(U)=U_{12}U_{13},
\end{equation}
where $U_{ij}$ is the usual ``leg notation'' (see \cite[Page 13]{KV1} and \cite[Section 7.1.2]{Tim}).

Let $\crp(\BG)$ denote the collection of unitary corepresentations of $\BG$.

For $U,\ V\in \crp(\BG)$, $T$ is called an \emph{intertwiner} between $U$ and $V$, and we write $T\in\intw(U,V)$, if $T\in\CL(\KH_U,\KH_{V})$ such that $$T(\id\otimes\omega)(U)=(\id\otimes\omega)(V)T,\ \text{for any}\ \omega\in L^1_*(\BG).$$
We say that $U$ is \emph{unitarily equivalent} to $V$ and write $U\cong V$, if there exists $T\in\intw(U,V)$ such that $T$ is a unitary.

%\medskip
\subsection{Universal quantum group}\label{subsec:univ-quant-gp}

The universal quantum group $C^*$-algebra of $\BGd$ is denoted by $(C_0^\ru(\BGd),\widehat{\Delta}^\ru)$ (see \cite[Section 4 and 5]{Kus}).
As shown in \cite[Proposition 5.2]{Kus}, there exists a unitary $$V_\BG^\ru\in M(C_0^\ru(\widehat{\BG})\otimes C_0(\BG))$$ that implements a bijection between unitary corepresentations $U$ of $\BG$ on $\KH$ and non-degenerate $*$-representations $\pi_U$ of $C_0^\ru(\widehat{\BG})$ on $\KH$ through the correspondence
$$U = (\pi_U\otimes\id)(V_\BG^\ru).$$
The identity $\Bo_\BG=1\otimes \Bo$ of $\CL(\BC)\otimes M(C_0(\BG))\cong M(\CK(\BC)\otimes C_0(\BG))\cong M(C_0(\BG))$ is a trivial unitary corepresentation of $\BG$ on $\BC$ and $\pi_{\Bo_\BG}$ is a character of $C_0^\ru(\widehat{\BG})$.

As in the literature, we write $U\prec W$ when $\pi_U\prec \pi_W$ (see, e.g., \cite[Section 5]{BT} and Section \ref{subsec:ntn}).

%\medskip

\subsection{Contra-corepresentation}

Let $U$ be a unitary coreprsentation of $\BG$ on a Hilbert space $\KH_U$.
As in \cite[Page 871]{BT}, we define the \emph{contra-corepresentation} $\overline U$ of $U$ by
$$\overline{U}:=(\tau\otimes R)(U),$$
where $\tau$ is the \emph{canonical anti-isomorphism} from $\CL(\KH_U)$ to $\CL(\overline{\KH}_U)$ (with $\overline \KH_U$ being the conjugate Hilbert space of $\KH_U$) and $R$ is the unitary antipode on $C_0(\BG)$.
Then $\overline U$ is a unitary corepresentation of $\BG$ on $\overline \KH_U$.
Notice that, it is unique up to equivalence $\cong$,
and that $$\overline{\overline{U}}\cong U.$$

%\smallskip

If $W$ is another unitary corepresentation of $\BG$ on a Hilbert space $\KK$, we denote by $U\circledT W$ the unitary corepresentation
$U_{13}W_{23}$ on $\KH\otimes\KK$ and call it the \emph{tensor product} of $U$ and $W$.
In this case,
\begin{equation}\label{eqt:ten-prod}
\pi_{U\circledT W} = (\pi_{U}\otimes\pi_{W})\circ \widehat{\Delta}^\ru.
\end{equation}

%\bigskip

\section{Amenability, co-amenability and Bekka-type amenabilities}

%In this section, we mainly generalize Proposition \ref{prop:ng-amen-bekkaamen} (or see the original version \cite[Proposition 3.6]{Ng}) to the general case, and give some discussions on the Bekka amenability.
Let us first recall the following definitions of amenability and co-amenability of a locally compact quantum group.

\begin{defn}\label{def:coamen-amen}{\rm(\cite[Definition 3.1 and 3.2]{BT})} Let $\BG$ be a locally compact quantum group.

\smnoind
(a) We say that $\BG$ is \emph{co-amenable} if there exists a state $\epsilon$ of $C_0(\BG)$) such that $(\id\otimes\epsilon)\Delta=\id$.

\smnoind
(b) A \emph{left invariant mean} for a locally compact quantum group $\BG$ is a state $m$ on $L^{\infty}(\BG)$ such that $m(\omega\overline{\otimes}\id)\Delta=\omega(\Bo)m$, for all $\omega\in L^1(\BG)$. 
We say that a locally compact quantum group $\BG$ is \emph{amenable} if it has a left invariant mean.
\end{defn}

\begin{rem}\label{rem:right-inv-mean}
Similarly, we also can define:
a \emph{right invariant mean} for $\BG$ is a state $m$ on $L^{\infty}(\BG)$ such that $m(\id\overline{\otimes}\omega)\Delta=\omega(\Bo)m$, for all $\omega\in L^1(\BG)$.
Clearly, $m$ is a right invariant mean if and only if $m\circ R$ is a left invariant mean. Thus, $\BG$ is amenable if and only if it has a right invariant mean.
\end{rem}

Co-amenability may be characterized by the following equivalent formulations, which were obtained by E.\ B$\acute{e}$dos and L.\ Tuset in \cite{BT}.

\begin{thm}\label{thm:coamen}{\rm(\cite[Theorem 3.1]{BT})}
For a locally compact quantum group $\BG$, the following statements are equivalent:

\smnoind
(a) $\BG$ is co-amenable.

\smnoind
(b) The canonical surjective homomorphism $\Lambda: C_0^\ru(\BG)\rightarrow C_0(\BG)$ is an isomorphism.

\smnoind
(c) There exists a $^*$-character on the $C^*$-algebra $C_0(\BG)$.

\smnoind
(d) There exists a net of unit vectors  $\{\xi_i\}$ in $L^2(\BG)$ such that $$\lim_i \norm{W_{\BG}(\xi_i\otimes v)-(\xi_i\otimes v)}=0,\ \ \forall v\in L^2(\BG).$$
\end{thm}

\begin{rem}\label{rem:amen-ng}
Comparing Theorem \ref{thm:coamen} with \cite[Theorem 2.3]{Ng}, we easily see that the ``amenability" for a Kac algebra in the Ng's paper \cite{Ng00} is actually the co-amenability for its dual in the sense of Definition \ref{def:coamen-amen}.
\end{rem}

%\smallskip

Now, we extend Bekka amenability and weakly Bekka amenability introduced in \cite{Ng} to the general case.

\begin{defn}\label{defn:wcp-bekkaamen} For any $U\in \crp(\BG)$, we say that

\smnoind
(a) $U$ has \emph{WCP (weak containment property)} if $\Bo_\BG\prec U$ (equivalently, $\pi_{\Bo_\BG}\prec \pi_U$, see Section \ref{subsec:ntn} and \ref{subsec:univ-quant-gp}). The WCP is actually the property $(A)$ introduced in \cite[Proposition and Definition 2.4]{Ng}.

\smnoind
(b) $U$ is \emph{Bekka amenable} if $\pi_{\Bo_\BG}\prec \pi_{U\circledT \overline{U}}$ (equivalently, $\Bo_\BG\prec U\circledT \overline{U}$), i.e., $U\circledT \overline{U}$ has the WCP.

%\smnoind
%(c) $U$ is \emph{conjugated Bekka amenable} if $\pi_{\Bo_\BG}\prec \pi_{\overline{U}\circledT U}$ (equivalently, $\Bo_\BG\prec \overline{U}\circledT U$), i.e., $\overline{U}\circledT U$ has the WCP.

\smnoind
(c) $U$ is \emph{weakly Bekka amenable} if there exists a positive functional $M$ on $\CL(\KH_U)$ with $M(\id_{\KH_U})=1$ such that
$$M[(\id_{\KH_U}\otimes\omega)(\alpha_U(T))]=M(T),$$ for any positive functional $\omega\in L^1(\BG)$ with $\omega(\Bo)=1$ and $T\in\CL(\KH_U)$, where $$\alpha_U(T):=U(T\otimes\Bo)U^*$$ is called a \emph{coaction} of $\BG$ on $\CL(\KH_U)$. Those $M$ satisfying the above condition are called \emph{$\alpha_U$-invariant means}.
\end{defn}

\begin{rem}\label{rem:conjbekkaamen} (a) Let $\BG$ be a locally compact quantum group of Kac type, and $U$ be an arbitrary finite dimensional unitary corepresentation of $\BG$.  
In \cite[Proposition 3.10]{Ng}, C.-K. Ng proved that $U$ is Bekka amenable.  
Clearly, so is $\overline{U}$, since $\overline{U}$ is also finite dimensional.

\smnoind
(b) When $\BG$ is actually a locally compact group $G$, its reduced $C^*$-algebraic quantum group $C_0(G)$ is commutative. 
It can be obviously seen from the commutativity that, for any two $U,\ V\in\crp(\BG)$, we have $U\circledT V=V\circledT U$. So, we have $U\circledT\overline{U}=\overline{U}\circledT U$, which implies that $U$ is Bekka amenable if and only if $\overline{U}$ is Bekka amenable. In this case, Bekka amenability is in fact the amenability for unitary representations in Bekka's paper \cite{Bek}.
\end{rem}

%\begin{rem}\label{rem:bekkaamen-conjbekkaamen-wcp}
%(a) For any two elements $U,\ U'\in\crp(\BG)$ and $U\cong U'$, it is obvious that if $U$ has the WCP, then $U'$ also has the WCP.

%\smnoind
%(b) For any $U\in\crp(\BG)$, it is easily known that $U$ is Bekka amenable if and only if $U\circledT\overline{U}$ has the WCP, and $U$ is conjugate Bekka amenable if and only if $\overline{U}\circledT U$ has the WCP.
%\end{rem}

\begin{thm}\label{thm:wcp}{\rm (\cite[Theorem 5.2]{BT} and \cite[Proposition and Definition 2.4]{Ng})}
Let $\BG$ be a locally compact quantum group and consider $U\in\crp(\BGd)$.
Then the following are equivalent:

\smnoind
(a) $U$ has the WCP.

\smnoind
(b) There exists a state $\psi$ on $\CL(\KH_U)$ such that $\psi(\id\overline{\otimes}\omega)(U)=\omega(\Bo)$, for $\omega\in L^1(\BGd)$.

\smnoind
(c) There exists a net $\{\xi_i\}$ of unit vectors in $\KH_U$ such that $$\lim_i \norm{U(\xi_i\otimes v)-(\xi_i\otimes v)}=0,\ \text{for all } v\in L^2(\BG).$$
\end{thm}

\begin{cor}\label{cor:coamen-wcp}
A locally compact quantum group $\BG$ is co-amenable if and only if $W_\BG$ has the WCP as a unitary corepresentation of $\BGd$.
\end{cor}

\begin{prf}
Since $W_\BG$ can be viewed as an element in $\crp(\BGd)$, the corollary easily follows from Theorem~\ref{thm:coamen}(d) and Theorem~\ref{thm:wcp}(c).
\end{prf}

The WCP is stable under some operations, for example, contra-gredient and tensor product.

%\smallskip

For any $U\in\crp(\BG)$, we denote by $\Bo_U$  the trivial unitary corepresentation $\id_{\KH_U}\otimes \Bo$ of $\BG$ on $\KH_U$.

\begin{prop}\label{prop:wcp}{\rm (\cite[Proposition 5.3]{BT})}
Suppose that $\BG$ is a locally compact quantum group and consider $U,\ V\in\crp(\BG)$.

\smnoind
(a) If $U$ has the WCP, then so does $\overline{U}$.

\smnoind
(b) If both of $U$ and $V$ have the WCP, then so does $U\circledT V$.

\smnoind
(c) If $U\circledT \Bo_{V}$ or $\Bo_{V}\circledT U$ has the WCP, then so does $U$.
\end{prop}

Next, we present a lemma and a proposition.
These results are probably known. Since we have not found them or their proof explicitly stated in the literature, we give complete arguments for the benefit of the reader.

\begin{lem}\label{lem:contra-tensor}
Let $U$ and $V$ be two unitary corepresentations of $\BG$. Then one has
$$\overline{U\circledT V}\cong\overline V\circledT\overline U.$$
\end{lem}

\begin{prf}
Since the unitary antipode $R$ is a $^*$-anti automorphism, one has
\begin{eqnarray*}
\overline{U\circledT V}
%&=&(\tau\otimes\tau\otimes R)(U_{13}V_{23})=(\id\otimes\tau\otimes R)(V_{23})(\tau\otimes\id\otimes R)(U_{13})=\overline{V}_{23}\overline{U}_{13}  \\
&=&(\tau\otimes\tau\otimes R)(U_{13}V_{23})\\
&=&(\tau\otimes R)(V)_{23}(\tau\otimes R)(U)_{13}\\
&=&\overline{V}_{23}\overline{U}_{13}=(\Sigma_{12}\otimes\Bo)\overline{V}_{13}\overline{U}_{23}(\Sigma_{12}\otimes\Bo)\\
&=& (\Sigma_{12}\otimes\Bo)(\overline{V}\circledT \overline{U})(\Sigma_{12}\otimes\Bo).
\end{eqnarray*}
So, it follows that, for any $\omega\in L^1_*(\BG)$, we have that
$$\Sigma(\id\otimes\omega)(\overline{U\circledT V})= (\id\otimes\omega)(\overline{V}\circledT \overline{U})\Sigma,$$
since $\Sigma^2=\id$. This implies that the unitary $\Sigma$ lies in $\intw(\overline{U\circledT V},\overline V\circledT\overline U)$. Hence, the Lemma holds.
\end{prf}

The following proposition is usually called \emph{the absorption principle}, which is the generalization of Fell's absorption principle for locally compact groups. E. B$\acute{e}$dos, R. Conti and L. Tuset in their paper \cite{BCT} proved the analogue in algebraic quantum groups (see \cite[Proposition 3.4]{BCT}).

\begin{prop}\label{prop:absorb}
Let $\BG$ be a locally compact quantum group. For any $U\in\crp(\BGd)$, one has that $$U\circledT W_{\BG}\cong\Bo_U\circledT W_\BG.$$
\end{prop}

\begin{prf}
Let $U$ be an arbitrary unitary corepresentation of $\BGd$. Set $T$ to be the image of $U$ on $\CL(\KH_U\otimes L^2(\BG))$.

For any $\omega\in L^1_*(\BGd)$, one has
\begin{eqnarray*}
T(\id\otimes\omega)(U\circledT W_{\BG})
& = & (\id\otimes\id\otimes\omega)(U_{12}U_{13}(W_{\BG})_{23})\\
& = & (\id\otimes\id\otimes\omega)((W_{\BG})_{23}U_{12})\\
& = & (\id\otimes\id\otimes\omega)((\id_{\KH_U}\otimes 1)_{13}(W_{\BG})_{23}U_{12})\\
& = & (\id\otimes\omega)(\Bo_U\circledT W_\BG)T,
\end{eqnarray*}
where the second ``$=$" comes from the pentagonal relation (see \cite[Definition A.1]{BS}): $U_{12}U_{13}(W_{\BG})_{23}=(W_{\BG})_{23}U_{12}$ for any $U\in\crp(\BG)$.

The above calculation implies that $T$ is a unitary intertwiner between $U\circledT W_{\BG}$ and $\Bo_U\circledT W_\BG$. Hence, the equivalence of $U\circledT W_{\BG}$ and $\Bo_U\circledT W_\BG$, as two elements in $\crp(\BGd)$, is obtained.
\end{prf}

\begin{cor}\label{cor:conj-absorb}
Let $\BG$ be a locally compact quantum group. For any $U\in\crp(\BGd)$, one has that
$$\overline{W}_{\BG}\circledT U\cong\overline{W}_{\BG}\circledT\Bo_U.$$
\end{cor}

\begin{prf}
For any $U\in\crp(\BGd)$, it is obvious that $\overline{U}$ is also in $\crp(\BGd)$.
By proposition \ref{prop:absorb}, one has that $$\overline{U}\circledT W_{\BG}\cong\Bo_{\overline{U}}\circledT W_\BG.$$ Hence, since $\Bo_{\overline{U}}\cong\overline{\Bo_U}$, by Lemma \ref{lem:contra-tensor}, we have that
$$\overline{W}_{\BG}\circledT U\cong \overline{\overline{U}\circledT W_\BG}\cong\overline{\Bo_{\overline{U}}\circledT W_\BG}\cong\overline{\overline{\Bo_U}\circledT W_\BG}\cong \overline{W}_{\BG}\circledT\Bo_U.$$
\end{prf}

\begin{cor}\label{cor:bekkaamen-conjbekkaamen}
Let $\BG$ be a locally compact quantum group. If $\overline{W}_\BG$ is Bekka amenable as a unitary corepresentation of $\BGd$, then $W_\BG$ is also Bekka amenable. 
%In other words, if $\overline{W}_\BG$ is Bekka amenable as a unitary corepresentation of $\BGd$, then $\overline{W}_\BG$ is conjugated Bekka amenable.
\end{cor}

\begin{prf}
Consider $W_\BG$ as a unitary corepresentation of $\BGd$. If $\overline{W}_\BG$ is Bekka amenable, then $\overline{W}_\BG\circledT W_\BG$ has WCP.
Hence, combining Corollary \ref{cor:conj-absorb} and Proposition~\ref{prop:wcp}(c), we have that $\overline{W}_\BG$ has WCP as a unitary corepresentation of $\BGd$. Hence, by Proposition \ref{prop:wcp}(a) and (b), both $W_\BG$ and $W_\BG\circledT \overline{W}_\BG$ also have the WCP, i.e. $W_\BG$ is Bekka amenable.
\end{prf}

Using the results above and the concept of the WCP, we can get a characterization of co-amenability for locally compact quantum groups.

\begin{prop}\label{prop:coamen-conjbekkaamen}
Let $\BG$ be a locally compact quantum group. The following statements are equivalent:

\smnoind
(a) $\BG$ is co-amenable.

%\smnoind
%(b) $W_\BG$ is conjugated Bekka amenable as a unitary corepresentation of $\BGd$.

\smnoind
(b) $\overline{W}_\BG$ is Bekka amenable as a unitary corepresentation of $\BGd$.
\end{prop}

\begin{prf}
First, assume that $\BG$ is co-amenable, that is, $W_{\BG}$ has the WCP by Corollary \ref{cor:coamen-wcp}.
Using the assertions (a) and (b) of Proposition~\ref{prop:wcp}, we know that $\overline{W}_\BG$ has the WCP and so does $\overline{W}_\BG\circledT W_\BG$.
Hence, by Definition~\ref{defn:wcp-bekkaamen}(b), $\overline{W}_{\BG}$ is Bekka amenable as a unitary corepresentation of $\BGd$.

Conversely, if $\overline{W}_{\BG}$ is Bekka amenable, then $\overline{W}_\BG\circledT W_\BG$ has the WCP. Considering $\overline{W}_\BG$ as a unitary corepresentation of $\BGd$, it follows from Proposition~\ref{prop:absorb} that $\Bo_{\overline{W}_\BG}\circledT W_\BG$ has the WCP.
Consequently, by Proposition~\ref{prop:wcp}(c), we know that $W_\BG$ has the WCP. Therefore, using Corollary~\ref{cor:coamen-wcp} again, we know that $\BG$ is co-amenable.
\end{prf}

\begin{cor}\label{cor:coamen-bekkaamen}
Let $\BG$ be a locally compact quantum group. If $\BG$ is co-amenable, then $W_\BG$ is Bekka amenable as a unitary corepresentation of $\BGd$.
\end{cor}

\begin{prf}
It follows directly from Corollary \ref{cor:bekkaamen-conjbekkaamen} and Proposition \ref{prop:coamen-conjbekkaamen}.
\end{prf}

In \cite{Ng}, using Bekka amenability of the fundamental multiplicative unitary, C.-K. Ng gave a characterization of ``amenability" of a Kac algebra. Using our terminology, we rewrite Ng's proposition as follows.

\begin{prop}\label{prop:ng-amen-bekkaamen}{\rm (\cite[Proposition 3.6]{Ng})}
Let $\BG$ be a locally compact quantum group of Kac type. Then $\BG$ is co-amenable if and only if $W_\BG$ is Bekka amenable as a unitary corepresentation of $\BGd$.
\end{prop}

As a direct consequence of Proposition \ref{prop:coamen-conjbekkaamen} and Proposition \ref{prop:ng-amen-bekkaamen}, the following corollary implies that, in the Kac case, $W_\BG$ is Bekka amenable if and only if $\overline{W}_\BG$ is Bekka amenable.
Note that the equivalence of (a) and (b) in this corollary is in fact Proposition \ref{prop:ng-amen-bekkaamen} proved by C.-K. Ng in \cite[Proposition 3.6]{Ng}.
We list these statements here just for comparison with the other results.

\begin{cor}\label{cor:kac-conjbekkaamen-bekkaamen}
Let $\BG$ be a locally compact quantum group of Kac type. Consider $W_\BG$ and $\overline{W}_\BG$ as two unitary corepresentations of $\BGd$. The following statements are equivalent:

\smnoind
(a) $\BG$ is co-amenable.

\smnoind
(b) $W_\BG$ is Bekka amenable.

%\smnoind
%(c) $W_\BG$ is conjugated Bekka amenable.

\smnoind
(c) $\overline{W}_\BG$ is Bekka amenable.

%\smnoind
%(e) $\overline{W}_\BG$ is conjugated Bekka amenable.
\end{cor}

In the following, we focus on weak Bekka amenability of unitary corepresentations. Using this property, we give another characterization for amenability, and generalizes some results on weak Bekka amenability in Ng's paper (see \cite[Proposition 3.4]{Ng}). Some proofs of these results below follow from similar lines of argument as that of \cite[Proposition 3.4]{Ng}. For completeness, we present the argument here.

\begin{prop}\label{prop:amen-weakbekkaamen}
Let $G$ be a locally compact quantum group. The following statements are equivalent:

\smnoind
(a) $\BG$ is amenable.

\smnoind
(b) the fundamental multiplicative unitary $W_\BGd$ of its dual group is weakly Bekka amenable as an element in $\crp(\BG)$.

\smnoind
(c) every $U\in\crp(\BG)$ is weakly Bekka amenable.
\end{prop}

\begin{prf}
To obtain that (a) implies (c), we first note that, by Remark \ref{rem:right-inv-mean} and amenability of $\BG$, there exists a right invariant mean $m$ on $\BG$.  Let $U$ be an arbitrary unitary corepresentation of $\BG$. For any positive functional $\omega$ on $\CL(\KH_U)$ with $\omega(\id_{\KH_U})=1$, we can define a linear map $\Phi_\omega$ from $\CL(\KH_U)$ to $C_0(\BG)$ by $$\Phi_\omega(T)=(\omega\otimes\id)\alpha_U(T),\ \text{for any}\ T\in\CL(\KH_U).$$
Furthermore, one can easily show that $\Phi_\omega$ is a completely positive map such that $\Delta\circ\Phi_\omega=(\Phi_\omega\otimes\id)\circ\alpha_U$ and $\Phi_\omega(\id_{\KH_U})=1$. Thus, we have that $M=m\circ\Phi_\omega$ is an $\alpha_U$-invariant mean for $U$, and so $U$ is weakly Bekka amenable. By arbitrariness of  the choice of $U$ , the statement (c) holds.

It is clear that  (c) implies (b), since $W_\BGd$ can be viewed as a unitary corepresentation of $\BG$.

To show that (b) implies (a), assume that $W_\BGd$ is weakly Bekka amenable, and let $\omega$ be an $\alpha_{W_\BGd}$-invariant mean. Hence, statement (a) follows from the fact that the restriction $\omega|_{L^{\infty}(\BG)}$ is indeed a left invariant mean for $\BG$.
\end{prf}

As in the Kac case, Bekka amenability is still stronger than weak Bekka amenability in the general case.

\begin{prop}\label{prop:bekkaamen-weakbekkaamen}
Let $\BG$ be a locally compact quantum group and $U$ be any unitary corepresentation of $\BG$. If $U$ is Bekka amenable, then $U$ is weakly Bekka amenable.

%\smnoind
%(b) If $U$ is conjugated Bekka amenable, then $\overline{U}$ is weakly Bekka amenable.
\end{prop}

\begin{prf}
(a) If $U$ is Bekka amenable, then, we know that $U\circledT\overline{U}$ has the WCP. Consequently, by Theorem \ref{thm:wcp}(c), there exists a net of unit vectors $\{\xi_i\}\subset\KH_{U\circledT\overline{U}}$ such that, for any $v\in L^2(\BG)$,
$$\lim_i \norm{(U\circledT\overline{U})(\xi_i\otimes v)-\xi_i\otimes v}=\lim_i \norm{(U\circledT\overline{U})^*(\xi_i\otimes v)-\xi_i\otimes v}=0\ (*)$$
Then the net of the vector states $\{\omega_{\xi_i,\xi_i}\}$ has a subnet weak$^*$-convergent to some positive functional $m\in\CL(\KH_{U\circledT\overline{U}})^*$.

For any unit vector $v\in L^2(\BG)$ and $T\in\CL(\KH_{U\circledT\overline{U}})$, 
%since $$U\circledT\overline{U}\in M(\CL(\KH_{U\circledT\overline{U}})\otimes C_0(\BG))\subset\CL(\KH_{U\circledT\overline{U}}\otimes L^2(\BG)),$$ 
one has
\begin{eqnarray*}
&& m[(\id_{\KH_{U\circledT\overline{U}}}\otimes\omega_{v,v})(\alpha_{U\circledT\overline{U}}(T))]\\
&=& \lim_i\omega_{\xi_i,\xi_i}[(\id_{\KH_U}\otimes\omega_{v,v})(\alpha_{U\circledT\overline{U}}(T))]\\
&=& \lim_i\langle(T\otimes\id_{L^2(\BG)})(U\circledT\overline{U})^*(\xi_i\otimes v),(U\circledT\overline{U})^*(\xi_i\otimes v)\rangle\\
&=& \lim_i\langle(T\otimes\id_{L^2(\BG)})(\xi_i\otimes v),\xi_i\otimes v\rangle\quad(By\ Equation\,(*))\\
%&=& \lim_i\langle T(\xi_i),\xi_i\rangle\langle \id_{L^2(\BG)}(\xi),\xi\rangle\\
&=& \lim_i\omega_{\xi_i,\xi_i}(T)\norm{\xi}^2=m(T).
\end{eqnarray*}
Because every $\omega\in L^1(\BG)$ is a linear combination of $\omega_{v,v}$'s, the equalities above imply that $m$ is an $\alpha_{U\circledT\overline{U}}$-invariant mean.
Define the positive functional $M$ on $\CL(\KH_U)$ by $M(T)=m(T\otimes\id_{\KH_{\overline{U}}})$ for any $T\in\CL(\KH_U)$.
Then, we can obtain weak Bekka amenability of $U$ by checking that $M$ is indeed an $\alpha_U$-invariant mean as required.

%\smnoind
%(b) If $U$ is conjugated Bekka amenable, then, by Remark \ref{rem:conjbekkaamen}(a), $\overline{U}$ is Bekka amenable. So, we know that $\overline{U}$ is weakly Bekka amenable by the statement (a) above.
\end{prf}

%We end  with a simple application. By Proposition \ref{prop:coamen-conjbekkaamen}, Proposition \ref{prop:amen-weakbekkaamen} and Proposition \ref{prop:bekkaamen-weakbekkaamen}, we can obtain, in a new way, the following well-known result.

%\begin{cor}\label{cor:amen-coaamen} Let $\BG$ be a locally compact quantum group. If $\BGd$ is co-amenable, then $\BG$ is amenable. \end{cor}

%\begin{prf}
%If $\BGd$ is co-amenable, then $W_\BGd$ is conjugated Bekka amenable as an element in $\crp(\BG)$. Consequently, by Proposition \ref{prop:bekkaamen-weakbekkaamen}, $W_\BGd$ is weakly Bekka amenable. Hence, $\BG$ is amenable by Proposition \ref{prop:amen-weakbekkaamen}.
%\end{prf}

%\bigskip

Finally, we conclude, from the results above, that for any locally compact quantum group $\BG$, the following relation holds:
%\smnoind
\begin{quote}
co-amenability of $\BGd$ $\Leftrightarrow$ Bekka amenability of $\overline{W}_{\BGd}$ $\Rightarrow$ Bekka amenability of $W_{\BGd}$ $\Rightarrow$ weak Bekka amenability of $W_{\BGd}$ $\Leftrightarrow$ amenability of $\BG$,
\end{quote}
where ``$\Leftrightarrow$" means``equal to" and ``$\Rightarrow$" means ``imply".

\bigskip
\bigskip
%-----------------------------------------------------------------------------------------
{\bf Acknowledgments.} The author sincerely appreciates Professor Chi-Keung Ng (Nankai University) for valuable suggestions. Meanwhile, the author also thanks the referees for valuable comments which improve this paper a lot.
The author are supported by the Post-doctoral Scientific Research Grant of Shandong University (1090516300002).

\bibliographystyle{amsplain}

\end{document}